\documentclass[12pt]{amsart}%
\usepackage{amsfonts}
\usepackage{amsmath}
\usepackage{amssymb}
\usepackage{graphicx}%
\theoremstyle{definition}
\newtheorem*{acknowledgments}{Acknowledgments}
\renewcommand{\theenumi}{\alph{enumi}}
\renewcommand{\MR}[1]{\relax}
\begin{document}
\title{Kadison--Singer from mathematical physics: An introduction}
\author{Palle E. T. Jorgensen}
\thanks{This material is based upon work supported by the National Science Foundation under Grant No.\ 
DMS-0457581.}
\address{Department of Mathematics, The University of Iowa, 14 MacLean Hall, Iowa City,
IA 52242-1419, U.S.A.}
\email{jorgen@math.uiowa.edu}
\urladdr{http://www.math.uiowa.edu/\symbol{126}jorgen/}
\begin{abstract}
                    We give an informal overview of the Kadison--Singer
extension problem with emphasis on its initial connections to
\linebreak
Dirac's formulation of quantum mechanics.

Let $\mathcal{H}$ be an infinite dimensional separable Hilbert space, and $\mathcal{B}\left(  \mathcal{H}\right)  $ the
algebra of all bounded operators in $\mathcal{H}$. In the language of operator algebras,
the Kadison--Singer problem asks whether or not for a given MASA  $\mathcal{D}$ in $\mathcal{B}\left(  \mathcal{H}\right)  $,
every pure state on $\mathcal{D}$ has a unique extension to a pure state on $\mathcal{B}\left(  \mathcal{H}\right)  $. In
other words, are these pure-state extensions unique?

                    It was shown recently by Pete Casazza and co-workers
that this problem is closely connected to central open problems in other
parts of mathematics (harmonic analysis, combinatorics (via
Anderson pavings), Banach space theory, frame theory), and applications
(signal processing, internet coding, coding theory, and more).

\end{abstract}
\maketitle

\addtolength{\leftmargini}{-6pt}

\section{\label{Int}Introduction}

This is a contribution to the webpage for an AIM 2006 Workshop on the
Kadison--Singer problem. The posted text will become a permanent Introduction
for the record. The current version is written by Palle Jorgensen, following
lectures at the meeting by Dick Kadison. The Introduction stresses how the
mathematical context and the problem itself grew out of conceptual issues in
quantum mechanics.

While the Kadison--Singer problem from the original Kadison--Singer paper
\cite{KaSi59} arose from mathematical issues at the foundation of quantum
mechanics, it was found more recently to be closely connected to a number of
modern areas of research in mathematics and engineering.

Indeed $C^{\ast}$-algebra theory was motivated in part by the desire to make
precise fundamental and conceptual questions in quantum theory, e.g., the
uncertainty principle, measurement, determinacy, hidden variables, to mention
a few (see for example \cite{Em84}). The 1959 Kadison--Singer problem is and
remains a problem in $C^{\ast}$-algebras, and it has defied the best efforts
of some of the most talented mathematicians of our time. The AIM workshop was
motivated by recent discoveries where it was shown that the original problem
is equivalent to fundamental unsolved problems in a dozen areas of research in
pure mathematics, applied mathematics and engineering, including: operator
theory, Banach space theory, harmonic analysis, and signal processing. While
the other parts of the present website will discuss details on that, following
Kadison's presentation at the workshop, this introduction will spell out some
of the original motivation behind the problem at its conception.

This little Introduction to K--S is limited in scope. Here is what I tried to
do, and what I stayed away from.

A number of themes are only hinted at in passing, and they could easily be
expanded into a monograph. So all one can hope for is a list of pointers, and
some explanations at an intuitive level, like: \textquotedblleft what does our
mathematical definition of a \emph{state} have to do with Dirac's ideas from
physics?\textquotedblright\ This was covered in Kadison's presentation at the
meeting, but is now fleshed out a little. I tried hard to limit the length of
the Intro, and yet still make a little dent into the murkiness of ideas from
quantum mechanics. Feynman used to say: \textquotedblleft Anyone claiming to
understand quantum mechanics should be met with skepticism!\textquotedblright%
\ (Quoted from memory!)

I will only recall that in the mid 1920s, the period from 1925 to 27, the
pioneering papers of Heisenberg, of Schr\"odinger, and of Dirac shaped quantum
mechanics into the theoretical framework we now teach to students in physics
and mathematics; see the quote from Dirac at the end. Heisenberg's paper came
first (by a few weeks) and was based on the notion of transition
probabilities, transition between states which later took the form of ``rays''
in Hilbert space, or equivalently vector states. Via corresponding matrix
entries, from this emerged what became known as ``matrix mechanics.'' Only
Heisenberg didn't realize that his matrices were infinite. Hence later
additions by Max Born and John von Neumann introduced Hilbert space and
operator algebras in a systematic way that is now taken for granted.

Von Neumann's axioms agree well with Heisenberg's vision, but Schr\"odinger
formulated his equation as a partial differential equation (PDE), generalizing
the classical wave equation. This was in the context of function spaces,
$L^{2}$-spaces on a classical version of phase space, and it had the
appearance of being ``closer'' to ``classical'' views of physics.
Schr\"odinger's wave functions are elements in the $L^{2}$-spaces, hence
``wave mechanics.''

At first it was thought that the two proposed frameworks for quantum mechanics
were contradictory, one was ``right,'' but not the other! Fortunately von
Neumann quickly proved that the two versions are unitarily equivalent, and
since von Neumann's paper \cite{vNeu32b} and his book \cite{vNeu68} the
concept of unitary equivalence has played and continues to play a central role.

This little Introduction does not go into technical points regarding all the
more recent implications, connections and applications of the K--S idea:
frames, signals, etc. Others will do that; see however the Reference
Supplement at the end.

I aim at offering some intuition regarding ideas and terminology that
originate in quantum physics, and in von Neumann's response to Heisenberg,
Schrodinger, and Dirac. Most of this can be found in courses in functional
analysis and Hilbert-space theory, and operator algebras. The trouble is that
if all of this were to be done properly in the Intro, it could easily become a
ten-volume book set. All that is realistic is a brief little Invitation, and
even that isn't easy to do well.

I hope to clarify questions asked at the meeting. If the Intro bridges some of
the diverse fields represented at the meeting, that is a help. The
participants include a broad and diverse spectrum of fields from math and from
signal processing, but not too much math physics. Yet, a good part of the
motivation derives from Dirac's vision, and I feel that it makes sense to
accentuate this part of the picture.

Another reason something like this might help is that over the years such
central parts of operator theory as the spectral theorem and Heisenberg's
uncertainty principle have slipped out of the curriculum in departments of
both math and physics.

And yet the spectral theorem and some of Heisenberg's ideas are central to the
many diverse subjects touched by the K--S problem, certainly in harmonic
analysis and in signal processing.

\section{\label{MaPh}Math and physics}

While physics students learn of quantum mechanical \emph{states} with
reference to experiments in the laboratory, in functional analysis as it is
taught in mathematics departments, states are positive linear functionals
$\omega$ on a fixed $C^{\ast}$-algebra $\mathfrak{A}$. We will assume that
$\mathfrak{A}$ contains a unit, denoted $I$. In this context, the conditions
defining $\omega$ are $\omega\colon\mathfrak{A}\to\mathbb{C}$ linear,
$\omega\left(  I\right)  =1$ and $\omega\left(  a^{*}a\right)  \ge0$,
$a\in\mathfrak{A}$, which is the usual positivity notion in operator algebras.
As is well known, a linear functional $\omega$ on $\mathfrak{A}$ is a state if
and only if $\omega\left(  I\right)  =1$ and $\omega$ has norm $1$.

This characterization of states (as norm-$1$ elements $\omega$ in
the dual of $\mathfrak{A}$ with $\omega\left(  I\right)  =1$)
is a lovely little observation due to
Richard Arens in the early 1950s. It makes things so much easier: You can now
use the simplest Hahn--Banach theorem to produce all kinds of states for special purposes.

                      The following little picture (cited from \cite{Jor03})
illustrates with projective geometry the simplest instance of
Pauli spin matrices, and it offers a lovely visual version of the
distinction between pure states and mixed states.

Recall first that the
familiar two-sphere $S^{2}$
goes under the name \textquotedblleft the Bloch sphere\textquotedblright%
\ in physics circles (to Pauli, a point in $S^{2}$ represents the state of an
electron, or of some spin-$1/2$ particle, and the points in the open ball
inside $S^{2}$ represent mixed states), and 
points in $S^{2}$ are identified with
equivalence classes
of unit vectors in $\mathbb{C}^{2}$, where equivalence of vectors $\mathbf{u}$
and $\mathbf{v}$ is defined by $\mathbf{u}=c\mathbf{v}$ with $c\in\mathbb{C}$,
$\left\vert c\right\vert =1$. With this viewpoint, a one-dimensional
projection $p$ on $\mathbb{C}^{N}$ is identified with the equivalence class
defined from a basis vector, say $\mathbf{u}$, for the one-dimensional
subspace $p\left(  \mathbb{C}^{N}\right)  $ in $\mathbb{C}^{N}$. A nice
feature of the identifications, for $N=2$, is that if the unit-vectors
$\mathbf{u}$ are restricted to $\mathbb{R}^{2}$, sitting in $\mathbb{C}^{2}$
in the usual way, then the corresponding real submanifold in the Bloch sphere
$S^{2}$ is the great circle: the points $\left(  x,y,z\right)  \in S^{2}$
given by $y=0$. To Pauli, $S^{2}$, as it sits in $\mathbb{R}^{3}$, helps
clarify the issue of quantum obervables and states. Pauli works with three
spin-matrices for the three coordinate directions, $x$, $y$, and $z$. They
represent observables for a spin-$1/2$ particle. States are positive
functionals on observables, so Pauli gets a point in $\mathbb{R}^{3}$ as the
result of applying a particular state to the three matrices. The pure states
give values in $S^{2}$. Recall that pure states in quantum theory correspond
to rank-one projections, or to equivalence classes of unit vectors.

                In conclusion, for this little Pauli spin model, the pure
states are realized as points on the $2$-sphere $S^{2}$, while the mixed states
are points in the interior, i.e., in the open ball in $\mathbb{R}^{3}$, centered at zero,
and with radius $1$.

              We now turn to the question: What does the positivity part of
the mathematical definition of a ``state'' given above have to do with
physics and experiments?

             The answer lies in the way Heisenberg introduced probability
into quantum mechanics. We elaborate on this point in eqs.\ (\ref{eqONB.6}) and (\ref{eqONB.8}) below
in a special case.

              At an intuitive level, the states from operator algebras serve
to make precise an analogy to the more familiar distribution of a random
variable from classical probability theory: In this ``correspondence
principle'', random variables correspond to selfadjoint operators, that is,
operators $A$ ($= A^*$) generalize random variables to the non-commutative, or
operator algebraic, setup of quantum mechanics. This \emph{ansatz} is consistent,
since, by the spectral theorem, every selfadjoint operator is represented up
to unitary equivalence as multiplication by a real-valued measurable
function in a suitable $L^2$-space. 

              Hence, in quantum mechanics, observables are selfadjoint
elements in the particular $C^*$-algebra which is selected to model the system
to be studied. Now the positivity axiom: For each (mathematical) state
$\omega$ we are assigning a probability to measurements of observables
prepared in experiments, i.e., in associated experimental states $S$
(instruments, prisms, magnetic fields, etc.). The information contained in $S$
is condensed into $\omega$. Now for the probability distributions: Given an
interval $J$ on the real line, and given a physical state $S$, we must calculate
the probability of measuring a quantum-mechanical observable $A$ (e.g.,
position, momentum, etc.)\ attaining values in $J$ when it is measured in some
prescribed and prepared state $S$, or rather $\omega$. When $A$ is given, the
probabilities come from the spectral theorem applied to $A$: that is, in the
form of a direct integral decomposition as recalled in (\ref{eqONB.6}) below.

\section{\label{Sta}States \& representations}

The mathematical significance of \emph{states} lies in their relationship to
\emph{representations}. By a \emph{representation} of a $C^{\ast}$-algebra
$\mathfrak{A}$ we mean a homomorphism $\pi\colon\mathfrak{A}\rightarrow
\mathcal{B}\left(  \mathcal{H}\right)  $, i.e., $\pi$ is linear, $\pi\left(
ab\right)  =\pi\left(  a\right)  \pi\left(  b\right)  $, $\pi\left(  a\right)
^{\ast}=\pi\left(  a^{\ast}\right)  $, $a,b\in\mathfrak{A}$, where
$\mathcal{H}$ is some complex Hilbert space, and where $\mathcal{B}\left(
\mathcal{H}\right)  $ denotes the $C^{\ast}$-algebra of all bounded operators
on $\mathcal{H}$. Note that the space $\mathcal{H}$ depends on $\pi$. While an
abstract $C^{\ast}$-algebra has an inherent $C^{\ast}$-involution, $\ast$, the
involution on $\mathcal{B}\left(  \mathcal{H}\right)  $ is defined as the
\emph{adjoint}, i.e., $A\rightarrow A^{\ast}$ defined by%
\begin{equation}
\left\langle \,Au\mid v\,\right\rangle =\left\langle \,u\mid A^{\ast
}v\,\right\rangle ,\qquad u,v\in\mathcal{H}, \label{eqSta.1}%
\end{equation}
where, as is
customary in physics/quantum mechanics, $\left\langle \,\cdot\mid\cdot\,\right\rangle $ denotes the inner
product in $\mathcal{H}$.

For a given state $\omega$ on a $C^{\ast}$-algebra $\mathfrak{A}$, there is a
unique cyclic representation $\left(  \pi_{\omega},\mathcal{H}_{\omega}%
,\Omega\right)  $, where $\mathcal{H}_{\omega}$ is a Hilbert space, $\Omega
\in\mathcal{H}_{\omega}$, $\left\Vert \Omega\right\Vert =1$, and
\begin{equation}
\omega\left(  a\right)  =\left\langle \,\Omega\mid\pi_{\omega}\left(
a\right)  \Omega\,\right\rangle ,\qquad a\in\mathfrak{A}. \label{eqSta.2}%
\end{equation}
This is the so-called Gelfand--Naimark--Segal (GNS) representation. By taking
orthogonal direct sums, it follows that every $C^{\ast}$-algebra is faithfully
represented as a subalgebra of $\mathcal{B}\left(  \mathcal{H}\right)  $ for
some (\textquotedblleft global\textquotedblright) Hilbert space $\mathcal{H}$.

Some of the questions on the interface of math and physics relate to choices
of $C^{\ast}$-algebras which are right for quantum mechanics. To understand
this, recall that quantum-mechanical \emph{observables} (in the mathematical
language) are selfadjoint operators $A$ (i.e., $A=A^{\ast}$) in Hilbert space.
A choice of $C^{\ast}$-algebra $\mathfrak{A}$ implies a choice of observables%
\begin{equation}
\mathfrak{A}_{\operatorname*{sa}}:=\left\{  \,A\in\mathfrak{A}\mid A=A^{\ast
}\,\right\}  . \label{eqSta.3}%
\end{equation}

                   We mentioned the equivalence of Heisenberg's and Schr\"odinger's formulations, i.e.,
matrix mechanics and wave mechanics, took the axiomatic form of unitary equivalence via the Stone--von-%
Neumann uniqueness theorem. However this does not suffice for infinite systems. There is a second
version of equivalence which is tied more directly to $C^*$-algebras. It enters consideration in physics
when passing from a finite number of degrees of freedom to an infinite number. However, this extension
is not just a curiosity, and in fact is dictated by quantum statistical mechanics, and by quantum
field theory. For their axiomatic formulations, we refer to the books by David Ruelle \cite{Rue04}, and by
Alain Connes \cite{Con94}. Very briefly: in the infinite cases, it turns out that different particles are
governed by different statistics, e.g., bosons and fermions. In fact, the Stone--von-Neumann uniqueness
theorem is false for these infinite variants; false in the sense that the natural representations are
not unitarily equivalent. Specifically, a choice of statistics automatically selects an associated
$C^*$-algebra $\mathfrak{A}$ for the problem at hand, e.g., the $C^*$-algebra of the canonical commutation
relations (CCRs), or the canonical anticommutation relations (CARs).  In the infinite case, there are
issues about passing to the limit, from finite to infinite; but when the statistics and therefore the
$C^*$-algebra are chosen, then the relevant representations will typically not be unitarily equivalent.
Nonetheless, there are uniqueness theorems that take the form of $C^*$-isomorphisms. As it turns out, in
the infinite case, these $C^*$-isomorphisms are not unitarily implemented.

\section{\label{ONB}Orthonormal bases (ONB)}

If the chosen Hilbert space $\mathcal{H}$ is separable, we may index
orthonormal bases (ONBs) in $\mathcal{H}$ by the integers $\mathbb{Z}$, i.e.,
$\left\{  \,e_{n}\mid n\in\mathbb{Z}\,\right\}  $. Such a choice $\left\{
\,e_{n}\mid n\in\mathbb{Z}\,\right\}  $ of ONB fixes a subalgebra
$\mathcal{D}\subset\mathcal{B}\left(  \mathcal{H}\right)  $ of operators $A$
which are simultaneously diagonalized by $\left\{  e_{n}\right\}  $, i.e.,
\begin{equation}
Ae_{n}=\lambda_{n}e_{n},\qquad n\in\mathbb{Z},\;\lambda_{n}\in\mathbb{C},
\label{eqONB.4}%
\end{equation}
with the sequence $\left(  \lambda_{n}\right)  $ depending on $A$. Hence
$\mathcal{D}$ is a \textquotedblleft copy\textquotedblright\ of $\ell^{\infty
}\left(  \mathbb{Z}\right)  $.

Operators of the form (\ref{eqONB.4}) may be written as%
\begin{equation}
A=\sum_{n\in\mathbb{Z}}\lambda_{n}\left\vert e_{n}\right\rangle \left\langle
e_{n}\right\vert , \label{eqONB.5}%
\end{equation}
where we use Dirac's notation for the rank-$1$ projection $E_{n}$ with range
$\mathbb{C}\,e_{n}$.

Caution to mathematicians: Eq.\ (\ref{eqONB.5}) is physics lingo, a favorite
notation of Dirac (\cite{Dir47,Dir39} quoted in \cite{KaSi59}). Now let us
facilitate the translation from physics lingo to (what has now become) math
notation. The representation in Eq.\ (\ref{eqONB.5}) is how the physicist
P.A.M. Dirac thought of diagonalization. Since the ``bras'' and the ``kets''
may be confusing to mathematicians, we insert explanation.

The notation used here is called Dirac's bra-ket (inner product), or ket-bra
(rank-one operator) notation; and it is adopted in the physics community, and
used in physics books.

Abstract considerations of Hilbert space are facilitated by Dirac's elegant
bra-ket notation, which we shall adopt. It is a terminology which makes basis
considerations fit especially nicely into an operator-theoretic framework: If
$\mathcal{H}$ is a (complex) Hilbert space with vectors $x$, $y$, $z$, etc.,
then we denote the inner product as a Dirac bra-ket, thus $\left\langle \,x
\bigm| y\,\right\rangle \in\mathbb{C}$. In contrast, the rank-one operator
defined by the two vectors $x$, $y$ will be written as a ket-bra, thus $E =
\left\vert \vphantom{\bigm|}x\right\rangle \left\langle
\vphantom{\bigm|}y\right\vert $. Hence $E$ is the operator in $\mathcal{H}$
which sends $z$ into $\left\langle \,y \bigm| z\,\right\rangle x$.

The general version of the spectral theorem for selfadjoint operators $A$ in
$\mathcal{H}$ takes the following form:%
\begin{equation}
A=\int_{\mathbb{R}}\lambda\,E\left(  d\lambda\right)  , \label{eqONB.6}%
\end{equation}
where $E\left(  \,\cdot\,\right)  $ is a projection-valued measure defined on
the sigma-al\-ge\-bra of all Borel subsets $\mathcal{B}$ of $\mathbb{R}$.
Specifically, for each $S\in\mathcal{B}$,%
\begin{equation}
E\left(  S\right)  ^{\ast}=E\left(  S\right)  =E\left(  S\right)  ^{2}.
\label{eqONB.7}%
\end{equation}

If a vector $v\in\mathcal{H}$, $\left\Vert v\right\Vert =1$, represents a
state (in fact a \emph{pure} state on $\mathcal{B}\left(  \mathcal{H}\right)
$) then%
\begin{equation}
\mathcal{B}\ni S\longmapsto\left\langle \,v\mid E\left(  S\right)
v\,\right\rangle =\left\Vert E\left(  S\right)  v\right\Vert ^{2}
\label{eqONB.8}%
\end{equation}
represents the probability of achieving a \emph{measurement} of the observable
$A$ with values in $S$ when an experiment is prepared in the state $v$,
written $\left\vert v\right\rangle $ in Dirac's terminology. If, further, the
system, prepared in the state corresponding to $v$, is designed to produce,
with certainty, $\lambda$, one of the possible values that a measurement of
the observable $A$ can yield (i.e., if the probability is $1$ that a
measurement of $A$ in this state will yield $\lambda$---an idealized extreme),
then $v$ is an eigenvector for $A$ corresponding to the eigenvalue $\lambda$.
If the measurement of $A$ in the general state $\omega$ of $\mathcal{B}\left(
\mathcal{H}\right)  $ yields $\lambda$ with certainty, we say that $\omega$ is
\emph{definite} on $A$. The condition for $\omega$ to be definite on the
observable $A$ is that $\omega\left(  A^{2}\right)  =\omega\left(  A\right)
^{2}$ \cite{KaSi59}.

While the three formulas (\ref{eqONB.6})--(\ref{eqONB.8}) are innocent-looking
assertions from pure mathematics, they grew out of Paul Dirac's endeavors in
making precise and extending the early formulations of quantum mechanics that
emerged from Werner Heisenberg, Erwin Schr\"{o}dinger, and Max Born; and later
the math physics schools in G\"{o}ttingen and in Copenhagen.

Section \ref{Con} below elaborates on these physics connections a bit more. To
summarize, the \textquotedblleft dictionary\textquotedblright\ is as follows.

\begin{enumerate}
\item \label{Dict(a)} Observable, e.g., momentum, position, energy,
spin${}\rightarrow{}$Selfadjoint operator, say $A$ in Hilbert space.

\item \label{Dict(b)} State (in the mathematical formulation as a positive
functional, say $\omega$)${}\rightarrow{}$Design and preparation of an
experiment in a laboratory, magnets, mirrors, prisms, radiation, scattering, etc.

\item \label{Dict(c)} Measurement (involving in its mathematical formulation
the spectral theorem as given in (\ref{eqONB.6}))${}\rightarrow{}$Application
of instruments to the observable $A$ as it is prepared in the state $\omega$.
\end{enumerate}

Caution: Note that an observable is not a number; it is a selfadjoint
operator. Because of Heisenberg's uncertainty relation, even a quantum
measurement is typically not really a number. Rather, it is the recording of a
probability distribution of a definite observable $A$ which is measured in a
specified state. This is what Eq.\ (\ref{eqONB.8}) is saying in the language
of functional analysis and operator theory.

\section{\label{Pur}Pure states}

Let $\mathfrak{A}$ be a $C^{\ast}$-algebra, and denote by $\Delta\left(
\mathfrak{A}\right)  $ the set of all states of $\mathfrak{A}$. From
functional analysis we know that $\Delta\left(  \mathfrak{A}\right)  $ is a
weak*-compact subset of $\mathfrak{A}_{1}^{\ast}:={}$the unit ball in the
dual. By Krein--Milman, we know that $\Delta\left(  \mathfrak{A}\right)  $ is
the closed convex hull of its \emph{extreme points}. The extreme points in
$\Delta\left(  \mathfrak{A}\right)  $ are known as the \emph{pure states} of
$\mathfrak{A}$.

The Kadison--Singer question is whether or not every pure state on
$\mathcal{D}$ has a unique pure-state extension to $\mathcal{B}\left(
\mathcal{H}\right)  $. But note that by Krein--Milman, it is only the
uniqueness part of the problem that is unresolved.

            Experts believe that the problem/conjecture is likely to be
``negative'' in the sense that there are pure states on $\mathcal{D}$ with multiple and
distinct pure-state extensions to $\mathcal{B}\left(
\mathcal{H}\right)  $. Specifically, this would mean that there are pure
states $\omega_{1}\neq\omega_{2}$ on $\mathcal{B}\left(
\mathcal{H}\right)  $ extending the same pure state on
$\mathcal{D}$. In other words, starting with such a pair of pure states on $\mathcal{B}\left(
\mathcal{H}\right)  $, that the
common restriction to $\mathcal{D}$ will define the same pure state $\rho$ on $\mathcal{D}$. Recall
that purity for states on $\mathcal{D}$ is equivalent to the multiplicative rule,
$\rho\left(
AB\right)  =\rho\left(  A\right)  \rho\left(  B\right)  $ for all
$A,B\in\mathcal{D}$.

           We now recall that the question of what such possible bifurcation
states $\rho$ on $\mathcal{D}$ might possibly look like concerns special points in the
Stone--\v{C}ech compactification $\beta\left(  \mathbb{Z}\right)  $
of the integers $\mathbb{Z}$, specifically points
in the corona $:= \beta\left(  \mathbb{Z}\right)  \setminus \mathbb{Z}$. We discuss this briefly below; and there will
be much more detail in a separate chapter.

Such a negative solution, if it exists, appears to hint at a \textquotedblleft
strange\textquotedblright\ element of quantum-mechanical indeterminacy.

The problem is whether or not such a bifurcation may happen from some pure
states $\rho$ on $\mathcal{D}$.

In another of the presentations included elsewhere on the site for the
workshop, the Kadison-Singer problem is analyzed starting from the familiar
realization of the pure states on $\mathcal{D}$ as points in the
Stone--\v{C}ech compactification $\beta\left(  \mathbb{Z}\right)  $.

\section{\label{Con}Concluding remarks}

The notion of \textquotedblleft purity\textquotedblright\ for states has
significance in both physics and mathematics. In mathematics, pure states
enter into extremality considerations, in linear programming, variational
analysis, and in decomposition theory. Examples: (a)~Formula (\ref{eqONB.8}),
above, shows that numbers obtained in quantum measurements attain their
extreme values at pure states. (b)~The vectors $\left(  e_{n}\right)  $ in
formula (\ref{eqONB.5}) define pure states on $\mathcal{B}\left(
\mathcal{H}\right)  $, in fact the simplest kind of pure states, eigenvectors;
i.e., each vector $e_{n}$ from (\ref{eqONB.5}) defines the pure state
$\omega_{n}\left(  \,\cdot\,\right)  :=\left\langle \,e_{n}\mid{\cdot}%
\,e_{n}\,\right\rangle $ on $\mathcal{B}\left(  \mathcal{H}\right)  $, and for
the eigenvalues we have $\lambda_{n}=\omega_{n}\left(  A\right)  $. (c)~In
contrast, consideration of continuous spectrum and of Heisenberg's uncertainty
relations dictates the more elaborate formulas (\ref{eqONB.6})--(\ref{eqONB.8}%
) for the most general selfadjoint operators in Hilbert space.

Referring to formula (\ref{eqSta.2}) (Section \ref{Sta} above) from
mathematics, the counterpart of pure states in the GNS-correspondence between
states and representations (see, e.g., \cite{KaRi97b,Arv76}) is irreducibility
for the representation. Specifically: A given state $\omega$ of a $C^{\ast}%
$-algebra $\mathfrak{A}$ is pure if and only if the corresponding cyclic
representation $\pi_{\omega}$ of (\ref{eqSta.2}) is irreducible. And it is
known that irreducible representations in physics label elementary particles.
In thermodynamics, pure states label pure phases. More generally, physics is
concerned with composite systems and their decompositions into elementary
building blocks.

Since this picture involving the GNS-correspondence (\ref{eqSta.2}) includes
the case when the given $C^{\ast}$-algebra is the group $C^{\ast}$-algebra of
a locally compact group, convexity and direct integral theory yields an
abstract Plancherel formula \cite{Seg63} for the harmonic analysis of groups;
see \cite{Seg63} and \cite{Dix81}.

\section{\label{Edi}Editorial comment by Palle Jorgensen}

This is a draft of the Introduction to the Kadison--Singer IMA website
\texttt{http://www.aimath.org/WWN/kadisonsinger/}\,.
It grew out of a workshop at the AIM
institute (with NSF support) in Palo Alto in September, 2006. Part of the
workshop program is the creation of a permanent AIM website for the
Kadison--Singer Problem, and I was assigned to write the first draft of an
Introduction. Several things guided me:

\begin{enumerate}
\renewcommand{\theenumi}{\arabic{enumi}}

\item \label{guide(1)}Motivation and history. This means that I left out
mention of current trends, and that I did not include an updated Bibliography.
This is left to the other writers for the project.

\item \label{guide(2)}The exposition is close to the lecture presentation
Professor R.V. Kadison gave at the 2006 IMA workshop itself.

\item \label{guide(3)}It is close to the original K--S 1959 paper, and it
stresses Dirac's influence.

\item \label{guide(4)}Physics, Paul Dirac, Dirac's book, and Dirac's thinking
were central to the motivations. (The write-up should be understandable to
physicists, for example workers in quantum computation. It is clear that
Dirac's views and notation are popular in these circles!) These concerns mean
that the present K--S Introduction will not talk much about a lot of other
more recent applications of the K--S ideas, for example to signal processing.
Others will write about that.

\item \label{guide(5)}I wanted to bridge separate communities, pure
vs.\ applied, math vs.\ physics, etc. To do that, I say a few things in ways
that have become popular with physicists, but perhaps aren't especially
familiar to mathematicians.
\end{enumerate}

This current draft is a ``live'' document and is likely to undergo a few more iterations.

\begin{acknowledgments}
We thank the organizers Pete Casazza, Richard Kadison, and David Larson for
comments and especially for putting together a productive and enjoyable
workshop. And we are further grateful for enlightening comments from the
workshop participants, especially from Gestur \'{O}lafsson, Vern Paulsen, and
Gary Weiss. We further thank Brian Treadway for typesetting and for helpful suggestions.
\end{acknowledgments}

\begin{quotation}
I received an early copy of Heisenberg's first work a little before
publication and I studied it for a while and within a week or two I saw that
the noncommutation was really the dominant characteristic of Heisenberg's new
theory. It was really more important than Heisenberg's idea of building up the
theory in terms of quantities closely connected with experimental results. So
I was led to concentrate on the idea of noncommutation and to see how the
ordinary dynamics which people had been using until then should be modified to
include it.

\hfill---P. A. M. Dirac

\end{quotation}


\section*{References (supplement)}

An incomplete list of more recent papers regarding the Kadison--Singer
problem: Its implications in mathematics and its neighboring fields.

Comments:

(a)            Another link between physics and the more modern signal-%
pro\-ces\-sing/operator-theory versions of K--S exists via what is called
\emph{complementarity} in quantum mechanics. This is also a connection to some of
the other references to papers and books on $C^*$-algebras.

 Pairs of non-commuting operators are said to be in complementarity when the
partial information computed jointly from the pair is maximal compared to
the information content computed from the two individually. Example:
momentum and position. 

As a result, extensions of pure states on one MASA are needed because the
operators are non-commuting, and therefore do not have simultaneous spectral
resolutions. While the concept of complementarity dates back to Niels Bohr,
it has found more recent uses in harmonic analysis, quantum
information/computation, and in signal processing; see, e.g., \cite{BLMM06,%
DH06,MB06,PR04}, and more on the arXiv  \texttt{http://arxiv.org/}\,.

(b) Many and diverse papers make connections to K--S via pavings, via
combinatorics, matrix theory, logic, foundations, harmonic analysis, coding,
information theory, and via Banach-space theory. The idea is that pavings
and a lot of other parts of the big picture are or will be covered by other
authors.

(c)  The paper \cite{CCLV05} and others make the connection between K--S and the
Feichtinger conjecture. The Feichtinger conjecture asserts that every
bounded frame can be written as a finite union of Riesz basic sequences.
There are recent results on this for Weyl-Heisenberg frames; hence the
connection to complementarity.

(d) Eventually there will perhaps be a big combined bibliography, but this
is a continuing community project. For now, we make do with a minimal list.%
\vspace*{-21pt}%
\nocite{CaTr06,CFTW06,AkAn91,And80,And79a,And79b,Wea04,HaPa06}%

\bibliographystyle{amsplain}
\providecommand{\bysame}{\leavevmode\hbox to3em{\hrulefill}\thinspace}
\renewcommand{\refname}{\relax}

\end{document}